\documentclass[10pt]{article}

\usepackage[english]{babel}
\usepackage{amsmath,amssymb,epsfig,bbm}
\usepackage{dsfont}
\usepackage{psfrag}
\usepackage{framed}
\usepackage{color}
\usepackage{colortbl}
\usepackage[T1]{fontenc}

\usepackage{graphicx}

\definecolor{shadecolor}{gray}{0.925}
\definecolor{dunkelgrau}{rgb}{0.8,0.8,0.8}
\definecolor{hellgrau}{rgb}{0.95,0.95,0.95}

\begin{document}


\numberwithin{equation}{section}

\newtheorem{thmintro}{Theorem}
\newtheorem{thm}{Theorem}
\newtheorem{defi}{Definition}
\newtheorem{lem}{Lemma}
\newtheorem{prop}{Proposition}
\newtheorem{cor}{Corollary}
\newtheorem{aim}{Aim}

\newcommand{\pf}{\noindent{\bf Proof: }}
\newcommand{\1}{\mathds{1}}

\newcommand{\N}{\mathbbm{N}}
\newcommand{\Z}{\mathbbm{Z}}
\newcommand{\R}{\mathbbm{R}}
\newcommand{\Co}{\mathbbm{C}}
\newcommand{\pr}{\mathbbm{P}}
\newcommand{\G}{\mathcal{G}}
\newcommand{\B}{\mathcal{B}}
\newcommand{\Pot}{\mathcal{P}}
\newcommand{\F}{\mathcal{F}}

\newcommand{\Om}{\Omega}
\newcommand{\om}{\omega}

\newcommand{\ssubset}{\Subset}
\newcommand{\Int}{\text{int}}
\newcommand{\Ext}{\text{ext}}

\newcommand{\tinybullet}{\begin{tiny}\bullet\end{tiny} }
\newcommand{\smallbullet}{\begin{footnotesize}\bullet\end{footnotesize} }

\newcommand{\alt}[1]{\par\noindent \textbf{Alternative:} \textit{#1} \hfill \linebreak }
\newcommand{\bem}[1]{\marginpar{\textcolor{red}{\emph{#1}}}} 
\newcommand{\pC}{ \mathfrak{C} }
\newcommand{\W}[1]{ \mathfrak{W} \left( #1\right) }
\newcommand{\C}[1]{ \mathfrak{C} \left( #1\right) }
\newcommand{\CC}[2]{ \mathfrak{C}^{#1} \left( #2\right) }
\newcommand{\HC}[2]{ \mathfrak{HC}_{*#1} \left( #2\right) }
\newcommand{\HCC}[3]{ \mathfrak{HC}^{#1}_{*#2} \left( #3\right) }

\newcommand{\uparc}[1]{ \mathfrak{A}^{+*}_{\text{up}} (#1) }
\newcommand{\upcc}[1]{ \mathfrak{CC}^{+*}_{ #1 } }
\newcommand{\minupcc}[1]{ \mathfrak{CC}^{\min +*}_{{#1}^c} }
\newcommand{\maxupcc}[1]{ \mathfrak{CC}^{\max +*}_{{#1}^c} }

\newcommand{\ul}[1]{ {}^{\displaystyle\bullet} {#1}} 
\newcommand{\um}[1]{ \stackrel{\displaystyle\bullet}{{#1}} } 
\newcommand{\ur}[1]{ {#1}^{\displaystyle\bullet }  } 
\newcommand{\ml}[1]{\displaystyle\bullet\,\! {#1}} 
\newcommand{\mr}[1]{ #1 \!\displaystyle\bullet} 
\newcommand{\dl}[1]{{}_{\stackrel{}{ \displaystyle\bullet}} {#1}} 
\newcommand{\dm}[1]{ {}_{\displaystyle\stackrel{\displaystyle {#1}}{\displaystyle\bullet} }   }
\newcommand{\dr}[1]{ {#1}_{\stackrel{}{\displaystyle\bullet}}} 

\newcommand{\cw}[2]{ #1 \stackrel{1*}{\curvearrowright} #2 } 
\newcommand{\ccw}[2]{ #1 \stackrel{1*}{\rotatebox[origin=c]{180}{$\curvearrowleft$} } #2 }


\begin{center}
{\bf \Large 
Non-Coexistence of Infinite Clusters\\[1ex] in Two-Dimensional Dependent Site Percolation
}\\ 

\vspace{1 cm}
Sebastian Carstens\footnote[1]{
math@carstens.cc 
, Mathematisches Institut der LMU, Theresienstr. 39, 80333 Muenchen, Germany
}
\end{center}

\bigskip

\begin{abstract}
This paper presents three results on dependent site percolation on the square lattice. 
First, there exists no positively associated probability measure on $\{ 0,1\}^{\Z^2}$ with the following properties: 
a) a single infinite $0$cluster exists almost surely, 
b) at most one infinite $1*$cluster exists almost surely, 
c) some probabilities regarding $1*$clusters are bounded away from zero. 
Second, we show that coexistence of an infinite $1*$cluster and an infinite $0$cluster is almost surely impossible when the underlying probability measure is ergodic with respect to translations, positively associated, and satisfies the finite energy condition. 
The third result analyses the typical structure of infinite clusters of both types in the absence of positive association. 
Namely, under a slightly sharpened finite energy condition, the existence of infinitely many disjoint infinite self-avoiding $1*$paths follows from the existence of an infinite $1*$cluster. 
The same holds with respect to $0$paths and $0$clusters. 
\\


\end{abstract}

\bigskip

\section{Introduction}

In this article we consider interacting systems in which each site of the square lattice $\Z^2$ is equipped with a random ``spin'' taking value either $0$ or $1$. 
Two lattice sites are called adjacent if their Euclidean distance is $1$, and $*$adjacent if their distance is $1$ or $\sqrt{2}$. 
The lattice then splits into maximal connected or $*$connected subsets, called clusters resp.\! $*$clusters, on which the sites take the same spin. 
In this way we obtain clusters of $0$spins, called $0$clusters, and $*$clusters of $1$spins, called $1*$clusters. 
Likewise, one can speak of $0$paths and $1*$paths. 
A classical question in these models is ``Under which conditions on a probability measure does a single infinite $1*$cluster necessarily occupy so much space that none is left for another infinite cluster?''. 

This question was analysed by A. Gandolfi, M. Keane and L. Russo in \cite{GKR}. 
They showed that the existence of an infinite cluster of one type implies the finiteness of all clusters of the opposite type if the underlying probability measure on $\{ 0,1\}^{\Z^2}$ is ergodic\footnote{The assumption of ergodicity always refers to translations and includes the translation invariance.}, invariant under reflections in the coordinate axes, and positively associated in the following sense. 
\begin{defi}
We say a probability measure $\mu$ on $\{ 0,1\}^{\Z^2}$ is positively associated, if 
$$
\mu (A\cap B) 
\ge 
\mu (A)\mu(B) 
\, 
$$
for all increasing events $A$ and $B$.
An event $A$ is increasing, if $\xi\in A$ and $\eta\ge\xi$ (pointwise) implies $\eta\in A$.
\end{defi}
In particular, this implies the almost-sure uniqueness of an infinite cluster of a given type.
The latter problem was addressed in full generality by Burton and Keane \cite{BK}, who showed that merely translation invariance, together with the so-called finite energy condition, is sufficient for the existence of at most a single infinite cluster of a given type, either $0$ or $1*$. 
Note, however, that this result still allows the coexistence of an infinite $0$cluster and an infinite $1*$cluster. 
The finite energy condition, as discovered by Newman and Schulman in \cite{NS} and defined below, roughly says that every local configuration is compatible with anything that happens elsewhere.
\begin{defi}
A probability measure $\mu$ on $\{ 0,1\}^{\Z^2}$ satisfies the finite energy condition if, for every finite set $\Delta\subset\Z^2$, 
$$
\mu (\eta \text{ on } \Delta | \xi \text{ off } \Delta ) >0
$$
for all $\eta\in\{ 0,1\}^\Delta$ and $\mu $-a.e. $\xi\in\{ 0,1\}^{\Delta^c}$.
\end{defi}
Based upon the Burton-Keane uniqueness-theorem, Zhang simplified the proof of the result of Gandolfi et al. 
His elegant argument was first published by Grimmett in \cite{Gr} and also works  when reflection invariance is replaced by rotation invariance. 
For independent percolation, Bollobas and Riordan \cite{BR} adapted Zhang's idea to the case of k-fold symmetry for any $k\ge 2$. 

\begin{figure}[htb!]
\includegraphics[width=1\linewidth]{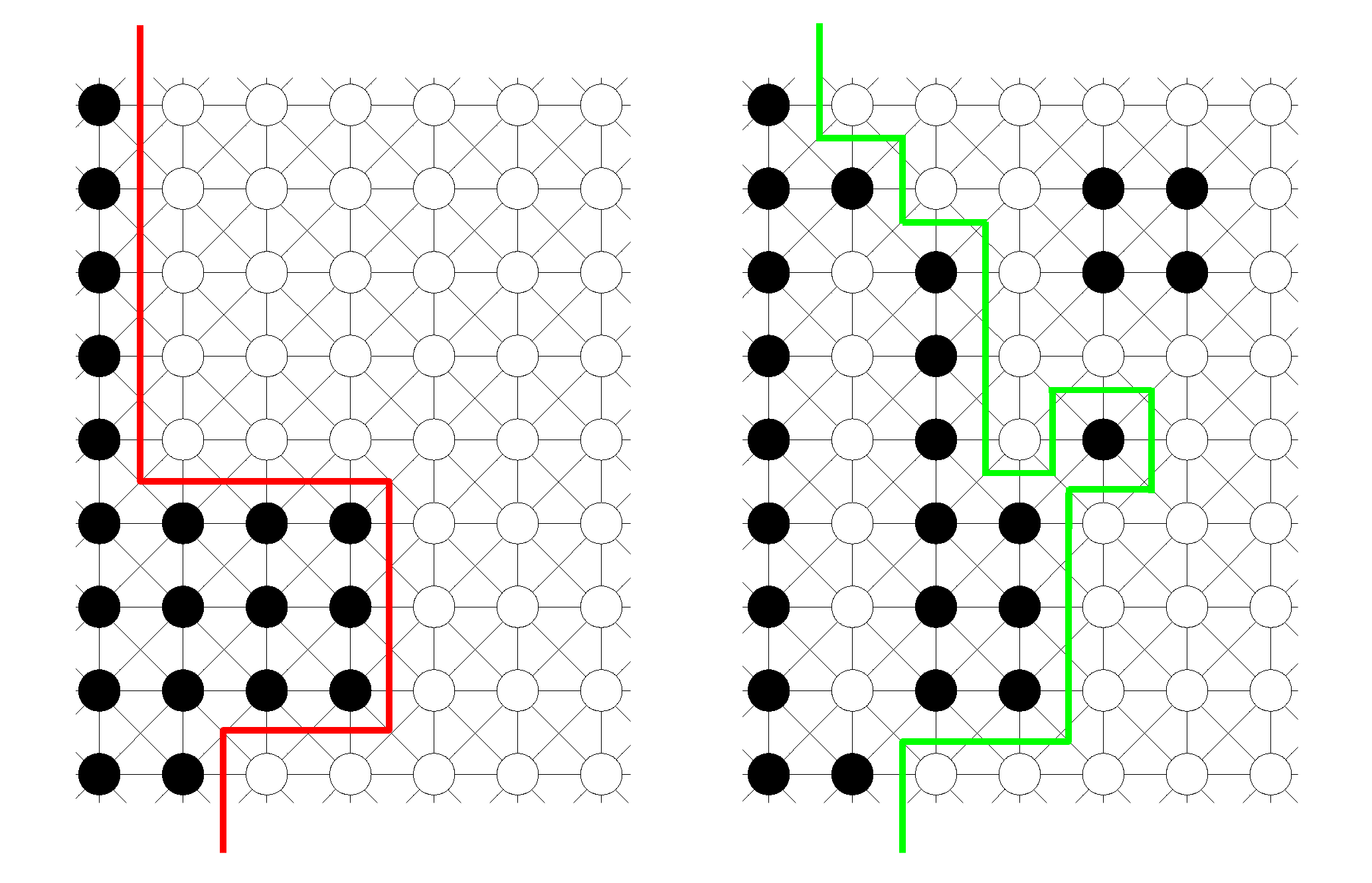} 
\caption{Black (resp.\! white) balls represent the nodes equipped with spins of value one (resp. zero). 
The horizontal, vertical, and diagonal lines from ball to ball represent the $*$edges. 
Sheffield's boundary is illustrated on the left side and the infinite boundary is illustrated on the right side.}
\label{Boundary}
\end{figure}

A rather different approach is due to Sheffield. 
Instead of assuming translation invariance and the finite energy condition, which are sufficient for the uniqueness of the infinite $0$clusters and the infinite $1$cluster, he directly assumed uniqueness and asked himself ``What conditions prevent the existence of an infinite $1$cluster if there already exists a single infinite $0$cluster'', see \cite[Theorem 9.3.1]{Sheff}.  
The formulation of his result requires the definition of the boundary between an infinite $0$cluster and an infinite $1$cluster. 
To this end, consider the event that the lattice splits up into a single infinite $0$cluster and a single infinite $1$cluster.
Then let each node of the infinite $1$cluster be the center of a unit square. 
Given all this, Sheffield defined the boundary between a single infinite $0$cluster and a single infinite $1$cluster, which is illustrated on the left side of figure \ref{Boundary}, as the topological boundary of the union of these squares. 
Using rather involved arguments, Sheffield was able to show that there exists no measure on $\{ 0,1\}^{\Z^2}$ such that 
\begin{itemize}
\item
it is positively associated; 

\item
with probability one, either all spins are $1$, or the lattice splits into a single infinite $0$cluster and a single infinite $1$cluster, and this second scenario happens with positive probability; 

\item
the distribution of the boundary between the two infinite clusters -- conditioned on their existence -- is translation-invariant.

\end{itemize}
In particular, almost surely relative to any translation invariant and positively associated measure on $\{ 0,1\}^{\Z^2}$, a single infinite $1$cluster and a single infinite $0$cluster cannot coexist; see \cite[Corollary 9.4.6]{Sheff}. 
His remarkable result dispenses completely with symmetries, but still relies on positive association as well as a kind of translation-invariance. 

We will analyse the same question as Sheffield,  considering infinite $1*$clusters in place of infinite $1$clusters, and obtain a theorem that dispenses with all kinds of  invariance under translations, reflections, or rotations. 
For this we need the following slight sharpening of the finite energy condition\footnote{Actually, we need much less. But this kind of condition is most convenient.}:
\begin{defi}
We say a probability measure $\mu$ on $\{ 0,1\}^{\Z^2}$ satisfies the bounded energy condition if for all $n\in\N$, there exists a strictly positive constant $c_n$ such that 
$$
\mu (\eta \text{ on } \Delta | \xi \text{ off } \Delta ) \ge c_n
$$
for all $\Delta\subset\Z^2$ with $|\Delta |=n$, all $\eta\in\{ 0,1\}^\Delta$, and for $\mu $-a.e. $\xi\in\{ 0,1\}^{\Delta^c}$.
\end{defi}
For example, the bounded energy condition holds for Gibbs measures relative to any shift-invariant and absolutely summable potential; cf.~\cite{Geo}. 

Now, we are ready to state our first result. 
\begin{thmintro}\label{mu}
There exists no probability measure $\mu $ on $\{0,1\}^{\Z^2}$ satisfying all of the following conditions: 
\begin{itemize}
\item[i)]
$\mu$ is positively associated;

\item[ii)]
$\mu$-almost surely there exists a single infinite $0$cluster; 

\item[iii)]
$\mu$-almost surely there exists at most one infinite $1*$cluster; 

\item[iv)]
$\mu$ satisfies the bounded energy condition; 

\item[v)]
the probability of the event ``a node $x$ belongs to the infinite $1*$cluster'' is bounded from below by a strictly positive constant not depending on $x\in\Z^2$.

\end{itemize}
\end{thmintro}

This theorem admits a version which is closer to Sheffield's theorem. 
Once again, the formulation requires the definition of the boundary between an infinite $0$cluster and an infinite $1*$cluster. 
To this end, consider the event that a single infinite $0$cluster and a single infinite $1*$cluster exist. 
Then fill the finite holes of the infinite $1*$cluster, i.e.\!, flip the spin of all $0$clusters $*$encircled by the infinite $1*$cluster. 
Now, let each node of this filled infinite $1*$cluster be the center of a square with sidelength $3/2$. 
Given all this, the infinite boundary, which is illustrated on the right side of figure \ref{Boundary}, is then defined as the topological boundary of the union of these squares. 
Note that, conditioned on the above coexistence event, the infinite boundary is always well-defined, since all nodes $*$adjacent to the infinite $1*$cluster are contained in the infinite $0$cluster.  
Furthermore, by definition, the infinite boundary indicates which side contains the infinite $0$cluster and differs from Sheffield's boundary between a single infinite $0$cluster and a single infinite $1$cluster. 
We usually interpret the infinite boundary as a curve. 

\begin{thmintro}\label{mu'}
There exists no probability measure $\mu' $ on $\{0,1\}^{\Z^2}$ which possesses all of the following properties: 
\begin{itemize}
\item[i)]
$\mu'$ is positively associated;

\item[ii)]
$\mu'$-almost surely there exists a single infinite $0$cluster; 

\item[iii)]
$\mu'$-almost surely there exists at most one infinite $1*$cluster; 

\item[iv')]
the distribution of the infinite boundary -- conditioned on its existence -- is translation-invariant; 

\item[v')]
with positive probability there exists an infinite $1*$cluster.

\end{itemize}
\end{thmintro}

The proof of these theorems requires only some fairly general properties of the matching pair $(\Z^2 ,\Z^{2*})$. 
As a consequence, we could replace $(\Z^2 ,\Z^{2*})$ by any matching pair $(M,M^* )$ with an underlying periodic mosaic $G$, where $M$ and $M^*$ are also periodic graphs. 
For the definition and a detailed discussion of matching pairs we refer to \cite{K}. 

The Burton-Keane uniqueness theorem, together with Theorem \ref{mu'}, implies our second result, namely the next corollary, which corresponds to the theorem of Gandolfi, Keane and Russo. 
In place of any kind of invariance under reflections or rotations, it takes advantage of the finite energy condition. 
\begin{cor}\label{rho}
Let $\rho $ be an ergodic and positively associated probability measure on $\{0,1\}^{\Z^2}$ satisfying the finite energy condition. 
Then the $\rho $-probability of coexistence of an infinite $1*$cluster and an infinite $0$cluster is zero. 
\end{cor}

As the assumption of positive association is often difficult to verify or does not hold at all, our third result analyses the structure of the event that an infinite $0$cluster or an infinite $1*$cluster exist under the sole condition of bounded energy. 
\begin{thmintro}\label{nu}
Let $\nu $ be a probability measure on $\{0,1\}^{\Z^2}$ satisfying the bounded energy condition. 
Then, $\nu $-almost surely on the event that an infinite $1*$cluster exists,
one can find infinitely many disjoint infinite self-avoiding $1*$paths. 
Analogously, infinitely many disjoint infinite self-avoiding $0$paths exist if an infinite $0$cluster exists. 
\end{thmintro}
If we assume -- in addition to the bounded energy -- the coexistence and uniqueness of the infinite $0$cluster and the infinite $1*$cluster, then Theorem \ref{nu} even yields the existence of infinitely many two-sided infinite self-avoiding $1*$paths. 
These two-sided infinite $1*$paths exhibit a natural order.
The first one $P_1$ is contained in the boundary of the infinite $0$cluster. 
The second of these two-sided infinite $1*$paths is contained in the boundary of the union of $P_1$, the infinite $0$cluster, and all finite $0$clusters adjacent to $P_1$. 
Since there exist infinitely many self-avoiding $1*$paths, this procedure can be continued indefinitely. 
So, the infinite $1*$cluster looks like wall bars. 
An analogous statement holds for the infinite $0$cluster and the lattice splitts into one $1*$wall bars and one $0$wall bars. 

Apart from that, Theorem \ref{nu} could also be useful as a first step towards a proof by contradiction of an analog of Corollary \ref{rho} which dispenses with positive association. 
If both the bounded energy condition and ergodicity hold, it seems to be counterintuitive that unique infinite clusters of both types coexist. 
For, on the one hand, the infinite $0$cluster is not allowed to intersect the intermediate space between the first and the $n$th two-sided infinite $1*$path as above, which has infinite ``length'' and ``width'' at least $n$.
On the other hand, ergodicity suggests that the infinite $0$cluster should be evenly spread over $\Z^2$. 

This intuition can be made rigorous under the further assumption of negative association, which means that any two increasing events are negatively correlated. 
Namely, subdivide the lattice into squares of the same size such that these squares can be interpreted as nodes of a new lattice.
Call two squares adjacent if their distance is one. 
Furthermore, call a square occupied if it is met by the infinite $0$cluster; otherwise it is called vacant. 
Given the coexistence, the size of the squares can be chosen so large that by exploiting the negative association, a standard path counting argument shows the finiteness of all vacant square-clusters. 
Let $N$ be a number that exceeds the diameter of the squares.
A two-sided infinite self-avoiding square-path is formed by the 
squares that are hit by the $N+1$th two-sided infinite self-avoiding $1*$path.
By choice of $N$, these squares are contained in the random set of nodes between the first and the $2N+1$th two-sided infinite self-avoiding $1*$path. 
Therefore, all of them are necessarily vacant, which is impossible because all clusters of vacant squares are finite. 

Once again, the proof of Theorem \ref{nu} only uses some fairly general properties of the matching pair $(\Z^2 ,\Z^{2*})$. 
So the proof also works for any matching pair $(M,M^* )$ with an underlying planar mosaic $G$ that satisfies the following two conditions:
The number of edges encircling a face is bounded, and the number of nodes adjacent to a node is bounded. 
Furthermore, the proof of the third result also works for higher dimensions.

\section{Some basic notations and definitions}
In this section we establish some basic notations needed throughout the paper. 
Our underlying set of nodes is always $\Z^2$. 

A node $x\in\Z^2$ is called adjacent to a set $B\subset\Z^2$ if $x\in\Z^2\setminus B$ and there exists a node $y\in B$ with $|x-y|=1$. 
Likewise, $x\in\Z^2$ is called $*$adjacent if $x\in\Z^2\setminus B$ and the Euclidean distance to $y\in B$ is 1 or $\sqrt{2}$. 
In particular, a node is not adjacent to itself. 
We define the boundary and $*$boundary of $B$ as $\partial B :=\{ x\in \Z^2 : x\text{ is } \text{adjacent to } B\} $ and $\partial^* B:=\{ x\in \Z^2 : x\text{ is } *\text{adjacent to } B\} $.

\begin{defi}\label{defi:path}
We call $(x_1,\ldots ,x_n)$, $n\ge 1$, a 
\begin{itemize}
\item 
path if 
$\,\forall 1\le i, j\le n: \, |i -j |=1 \Rightarrow x_i \text{ is adjacent to } x_j\, .$
\item 
self-avoiding path if it is a path and if 
$x_i=x_j\Rightarrow i=j$.
\end{itemize}
The node $x_1$ (resp.\! $x_n$) is called the starting (resp.\! ending) node. 
A sequence of nodes, $(x_i )_{i\ge 1}$, is an infinite self-avoiding path if, for all $n\in\N$, $(x_1,\ldots ,x_n)$ is a self-avoiding path. 
Further, a sequence of nodes, $(x_i )_{i\in\Z }$, is called a two-sided infinite self-avoiding path if the sequences $(x_i )_{i\ge 1}$ and $(x_i )_{i<1 }$ are two disjoint infinite self-avoiding paths, whose starting nodes are adjacent to each other.
\end{defi}
Accordingly, $*$paths and self-avoiding $*$paths are defined by replacing adjacent with $*$adjacent. 
A path hits $\Delta\subset\Z^2$ if one of its nodes belongs to $\Delta$. 

\begin{defi}\label{defi:circuit}
We call a path a 
\begin{itemize}
\item
circuit if its starting node is adjacent to its ending node.
\item
self-avoiding circuit if it is a circuit as well as a self-avoiding path. 
\end{itemize}
\end{defi}
Once again, we accordingly define $*$circuit and self-avoiding $*$circuit. 

By misuse of notation, a circuit or a path is often interpreted as a set. 
We write $\Delta\ssubset\Gamma$ to indicate that $\Delta$ is a finite subset of $\Gamma $. 

The interior of a circuit $C$, denoted $\Int C$, is the set of nodes in $\Z^2\setminus C$ that is enclosed by $C$. 
The exterior of $C$, $\Ext C$, is defined as $\Z^2 \setminus \left( C\cup\Int C\right)$. 
Whenever $\Delta\subset\Int C$ we say $C$ is a circuit around $\Delta$. 

\begin{defi}
Let $\sigma\in\{0,1\}^{\Z^2}$. 
We call a path $P$ a $0$path with respect to $\sigma$ if $P\subset\sigma^{-1} (0)$. 
\end{defi}
Self-avoiding $0$paths and $1*$paths and self-avoiding $1*$paths are analogously defined. 

Let $A,B,C\subset \Z^2$. 
We write $A\stackrel{0}{\longleftrightarrow} B$ in $C$ (resp.\! $A\stackrel{1*}{\longleftrightarrow} B$ in $C$) for the event that there exists a self-avoiding $0$path (resp.\! $1*$path) which belongs to $C$, starts in $A$ and ends in $B$. 
When $C=\Z^2$ the phrase ``in $\Z^2$'' is usually omitted. 
We exchange $B$ with $\infty$ to express that an infinite self-avoiding $0$path (resp.\! $1*$path), which is contained in $C$, exists and starts in $A$. 

%
\begin{defi}
Let $\sigma\in\{0,1\}^{\Z^2}$. 
A $0$cluster with respect to $\sigma $ is a subset $S\subset\sigma^{-1}(0)$ such that 
\begin{itemize}
\item[a)]
$\forall ~x,y\in S :~  x\stackrel{0}{\longleftrightarrow} y$ in $S$;

\item[b)]
$\nexists ~z\in S^c :~  z\stackrel{0}{\longleftrightarrow} S$, 

\end{itemize}
i.e.\!, $S$ is a maximal $0$connected subset of $\sigma^{-1}(0)$. 

\end{defi}

We define $1*$clusters and $0$circuits and $1*$circuits in the same way.

\section{Non-coexistence of infinite $1*$clusters and infinite $0$clusters}

In the main part of this section, we present a proof by contradiction of Theorem \ref{mu}.
Therefore, let $\mu $ be a probability measure on $\{ 0,1\}^{\Z^2}$ with the following properties: 
\begin{itemize}
\item[i)]
$\mu$ is positively associated;

\item[ii)]
$\mu$-almost surely there exists a single infinite $0$cluster; 

\item[iii)]
$\mu$-almost surely there exists at most one infinite $1*$cluster; 

\item[iv)]
$\mu$ satisfies the bounded energy condition; 

\item[v)]
there exists a constant $c>0$ such that $\mu (x\stackrel{1*}{\longleftrightarrow}\infty )\ge c$ for all $x\in\Z^2$.

\end{itemize}
Note that, nonetheless, an arbitrary number of finite clusters of both types could exist. 

We derive the contradiction in the following way: 
Let $\Delta\ssubset\Z^2$ be an arbitrary (but fixed) set containing the origin. 
We show that with probability at least $\epsilon >0$, the infinite $1*$cluster forms a $1*$circuit around $\Delta $. 
Moreover, $\epsilon$ does not depend on the choice of $\Delta$. 
So, if $\Delta\ssubset\Z^2 $ is large enough such that $\mu (\Delta\stackrel{0}{\longleftrightarrow} \infty ) \ge 1-\epsilon/2$, then the impossible event ``there exists a $1*$circuit around $\Delta$ as well as an infinite self-avoiding $0$path starts in $\Delta$'' has probability at least $\epsilon /2$, which is a contradiction.
Thus, an infinite $1*$cluster prevents the existence of an infinite $0$cluster. 

But how do we deduce the existence of $\epsilon $? 
Our strategy consists of the following three steps: 
First, if $x,y\in\Z^2$ are sufficiently far away from $\Delta$ the event that there exists a $1*$path from $x$ to $y$ in $\Delta^c$ occurs with probability at least $c^2/2$, where $c$ is as defined in property v). 
Second, a $1*$path from $x$ to $y$ in $\Delta^c$ could be either clockwise or counterclockwise coiled around the origin and the existence of both types implies the existence of a $1*$circuit around $\Delta$. 
Third, there exist $x,y\in\Delta^c$ such that with probability at least $c^2/4$, a clockwise $1*$path from $x$ to $y$ in $\Delta^c$ exists and with probability at least $c^2/4$ a counterclockwise $1*$path from $x$ to $y$ in $\Delta^c$ exists. 
This, together with the positive association and step two, implies that with probability at least $c^4/2^4 =\epsilon$, a $1*$circuit around $\Delta$ exists.

For the first step, we introduce a special self-avoiding $*$circuit, which consists of a $0$path and a $1*$path that are connected to form a $*$circuit.
\begin{defi}
Let $n,m\ge 0$ and $(x_1 ,\ldots ,x_n )$ be a self-avoiding $1*$path and $(y_1 ,\ldots ,y_m )$ be a self-avoiding $0$path with $x_n\stackrel{*}{\sim} y_1$ and $x_1\stackrel{*}{\sim} y_m$. 
We call the composition $(x_1 ,\ldots ,x_n, y_1 ,\ldots ,y_m )$ a (self-avoiding) mixed ${}^{1*}_{0}$circuit. 
\end{defi}
Note that a self-avoiding $1*$circuit or a self-avoiding $0$circuit is also a mixed ${}^{1*}_{0}$circuit. 

The purpose of this definition is the following: 
Let $\Delta\subset\Gamma\ssubset\Z^2 $ and $x,y\in\Gamma^c$. 
The existence of both a mixed ${}^{1*}_{0}$circuit in $\Gamma$ around $\Delta$ and a $1*$path from $x$ to $y$ implies that one can also find a $1*$path from $x$ to $y$ not hitting $\Delta$. 
Therefore, such a circuit ``shields'' $\Delta $ from ``outside'' $*$paths. 

\begin{lem}[Shield lemma]\label{+*0-hc}
For all $\Delta\ssubset\Z^2$, $\mu$-almost surely there exists a mixed ${}^{1*}_{0}$circuit around $\Delta$. 
\end{lem}

\pf 
It is sufficient to take $\Delta=\{ -d,\ldots ,d\}^2 $. 
We distinguish three cases.

First, we assume that all $1*$clusters meeting $\partial^* \Delta $ are finite. 
Then there exists a self-avoiding $0$circuit around $\Delta$, which in particular is a mixed ${}^{1*}_{0}$circuit in $\Delta^c$. 

The second case ``only finite $0$clusters meet $\partial^* \Delta $'' is solved analogously. 

Now, we turn our attention to the remaining case that the infinite $0$cluster and the infinite $1*$cluster (exist and) meet $\partial^*\Delta$. 
Thus, the infinite boundary (exists and) splits $\Z^2$ into two sets $S_0$ and $S_{1*}$ such that the one side $S_0$ consists of the infinite $0$cluster plus all its finite $*$holes, i.e.\!, $1*$clusters encircled by the infinite $0$cluster, and the other side $S_{1*}$ consists of the infinite $1*$cluster plus all its finite holes, i.e.\!, $0$clusters encircled by the infinite $1*$cluster. 

Because of the case assumption the infinite boundary hits $\partial^*\Delta$. 
Let $x,x'\in\partial^*\Delta \cap S_0 $ and $y,y'\in\partial^*\Delta \cap S_{1*} $ be the nodes such that the infinite boundary first enters $\partial^*\Delta$ between $x $ and $y $ and last exits $\partial^*\Delta$ between $x' $ and $y' $. 
In particular, $x$ is adjacent to $y$, $x'$ is adjacent to $y'$, the nodes $x, x'$ belong to the infinite $0$cluster, and $y, y'$ belong to the infinite $1*$cluster. 
Since all $1*$clusters in $S_0$ are finite and encircled by the infinite $0$cluster, which contains $x$ and $x'$, one can find a $0$path from $x$ to $x'$ in $S_0\cap \Delta^c $. 
Likewise, there exists a $1*$path from $y$ to $y'$ in $S_{1*}\cap \Delta^c $. 
The $0$path and the $1*$path are the necessary ingredients for a mixed ${}^{1*}_{0}$circuit around $\Delta$. 
Therefore, we have shown the existence in the third case. 

The lemma follows from the fact that one of these three cases almost surely occurs.
$\hfill\square$
\\

Note that only conditions ii) and iii) were used in this proof. 
The next lemma, which completes our first step towards proving Theorem \ref{mu}, relies on properties i) and v) in combination with the shield lemma. 
\begin{lem}\label{step2}
For all $\Delta\ssubset \Z^2$, there exists a set $\Gamma\ssubset\Z^2$ such that for all $x,y\in\Gamma^c $, the event ``$x$ and $y$ are $1*$connected in $\Delta^c$'' occurs with probability at least $c^2/2 $. 
\end{lem}

\pf 
Fix an arbitrary $\Delta\ssubset\Z^2$. 
Due to the shield lemma, we can choose $\Gamma\ssubset\Z^2$ such that with probability at least $1-c^2/2$, a mixed ${}^{1*}_{0}$circuit around $\Delta$ in $\Gamma $ exists. 
Let $x,y\in\Gamma^c$. 
The uniqueness of the infinite $1*$cluster yields the existence of a $1*$path from $x$ to $y$ as soon as $x$ and $y$ belong to this infinite $1*$cluster. 
Properties i) and v) imply that the latter event has probability at least $c^2$. 
Moreover, by the choice of $\Gamma$, we can conclude that with probability at least $c^2/2$, there exists in addition a mixed ${}^{1*}_{0}$circuit around $\Delta$ in $\Gamma $. 
Under these conditions, a $1*$path from $x$ to $y$ in $\Delta^c$ can be found. 
$\hfill\square$
\\


In our next step, the self-avoiding $*$paths from $x$ to $y$ off $\Delta\ssubset\Z^2$ are distinguished into two classes according to whether they run clockwise or counterclockwise around the origin. 
If $*$paths of both types exist, one can also find a $*$circuit around $\Delta$. 
To this end, we introduce the winding number around the origin, which for convenience will only be defined for polygons, i.e.\!, for piecewise linear continuous curves in $\R^2$. 

\begin{defi}
Let $\Delta$ be a finite set in $\R^2 $ that contains the origin and let $P : [0,1]\to \R^2\setminus\Delta $ be a polygon. 
We identify $\R^2$ with $\Co$ and rewrite $P(t)$ in polar form $P(t)=r(t)e^{i\theta (t) }$, where $\theta (.)$ is a continuous function. 
The winding number 
$$
\W{P}:= \frac{\theta (1)-\theta(0)}{2\pi }
$$
describes the fractional turns of the polygon around the origin. 
\end{defi}

We refer to \cite{B} for an alternative definition and elementary properties. 

Now, we are ready to define the two classes.
%
\begin{defi}
Let $x$ and $y$ be two nodes and let $P:[0,1]\to \R^2\setminus [-n,n]^2 $ be a polygon from $P(0)=x$ to $P(1)=y $. 
When $\W{P}$ is negative $P$ is called a clockwise polygon in $([-n,n]^2)^c $. 
When $\W{P}$ is positive $P$ is called a counterclockwise polygon in $([-n,n]^2)^c $. 
\end{defi}

The next lemma is a special case of the ``Topological Lemma'' in \cite{GKR}. 
It says that a $*$circuit exists if one can find a clockwise $*$path as well as a counterclockwise $*$path. 
Therefore, it concludes our second step.  
Obviously, $*$paths can be thought of as polygons.

\begin{lem}\label{cwccw}
Let $\Delta:=\{ -n,\ldots ,n\}^2$ and $x,y\in\Delta^c$. 
We assume that there exist a clockwise $*$path $P$ from $x$ to $y$ in $\Delta^c $ and a counterclockwise $*$path $Q$ from $x$ to $y$ in $\Delta^c$.
Then a $*$circuit around $\Delta$ in $P\cup Q$ exists.
\end{lem}

\pf
We consider the closed polygon $C(t):= P(2t)\1_{t<1/2}+ Q(2-2t)\1_{t\ge 1/2}$. 
Standard properties of the winding number yield 
$$
\W{C} 
= 
{\W{P}} - {\W{Q}}
\, ,
$$
which is negative, because the first summand is negative and the second one is positive. 
So the origin belongs to a bounded component of $\R^2\setminus C$. 
Consequently, there exists a $*$circuit around $\Delta $ that follows $P$ from $x$ to $y$ and then $Q$ backwards from $y$ to $x$. 
$\hfill\square$
\\


The aim of our third step is to verify the existence of two nodes $x,y$ such that the probabilities of the events ``there exists a clockwise $1*$path from $x$ to $y$ around $\Delta$'', in short $\cw{x}{y}$ around $\Delta$, and ``there exists a counterclockwise $1*$path from $x$ to $y$ around $\Delta$'', in short $ \ccw{x}{y} $ around $\Delta$, are bounded from below by a strictly positive constant, which does not depend on $\Delta$. 

The phrase ``$x$ is on the left side of $\Delta$'' means that one can find $d\in\N$ such that $x\in\{ (i,j): i\le -d\}$ and $\Delta\subset [-d,d]^2$ hold. 
Accordingly, ``a node is on the right side of $\Delta$'' is used. 

First we pursue the following idea: 
A $1*$path, that starts on the left side and ends on the right side of the origin, becomes a clockwise polygon when it is sufficiently shifted upwards. 

Note that the existential quantifier of the next lemma could be replaced with a universal quantifier, but stating the weaker version simplifies the modification for Theorem 2. 
\\

\begin{lem}\label{step3}
For all $\Gamma \ssubset\Z^2$, there exist a node $x$ on the left side and a node $y$ on the right side of $\Gamma$ such that 
\begin{align}
\exists h>0:~ &
\mu (\cw{x_h }{y_h } \text{ around }\Gamma  ) 
\ge  
c^2 /4
\label{cwh>0}
\\
\exists h<0:~&
\mu (\ccw{x_h }{y_h } \text{ around }\Gamma  ) 
\ge  
c^2 /4 
\label{cwh<0}
\, ,
\end{align}
where $x_h := x+(0,h)$ and $y_h := y+(0,h)$. 
\end{lem}

\pf
Since the proofs of \eqref{cwh>0} and \eqref{cwh<0} are obviously similar, we just verify \eqref{cwh>0}. 
The idea is more or less the same as in Lemma \ref{step2}.

Fix an arbitrary $\Gamma \ssubset\Z^2$ and choose $m\in\N$ such that $\Gamma\subset [-m,m]^2 $. 
Let $x:= (-m-1,0)$ and $y:= (m+1,0)$, which, therefore, are on the left respectively on the right side of $\Gamma$, and assume for contradiction that 
\begin{align}\label{contra1}
\forall h>0:\quad
\mu (\cw{x_h }{y_h } \text{ around }\Gamma  ) 
<c^2 /4
\, .
\end{align}
Let $P(h)$ be the shortest path from $x_h $ to $y_h $, i.e.\!, 
$$
P(h):= 
\left( (-m-1,h),(-m,h),\ldots, (m,h), (m+1,h) \right)
\, .
$$ 
The bounded energy condition ensures the existence of a constant $\delta >0$ such that with probability at least $\delta$, for all $h>0$, all spins of $P(h)$ take the value $1$. 
In particular, for all $h>m$
\begin{align}\label{h>m}
\mu \left( \cw{x_h }{y_h } \text{ around } \Gamma \right) 
\ge 
\delta
\, .
\end{align}
Let $\delta':=\delta c^2/ 4 $ and let $\Lambda\ssubset\Z^2$ be such that $\Lambda$ contains $\{-m,\ldots ,m \}^2$ and
\begin{align}
\label{0-Lambda}
\mu 
\left( 
\Lambda\stackrel{0}{\longleftrightarrow}\infty
\right)
>
1-\delta' /2
\, .
\end{align}
Due to Lemma \ref{step2}, there exists a square $\{-l,\ldots ,l\}^2 $ including $\Lambda $ such that for all $h>l$
$$
\mu 
\left(
x_h\stackrel{1*}{\longleftrightarrow}y_h \text{ around } \Lambda 
\right)
\ge 
c^2/2
\, .
$$
This, together with
$$
\{\cw{x_h}{y_h} \text{ around }\Lambda \}\cup \{\ccw{x_h}{y_h} \text{ around }\Lambda \}
=
\{ x_h\stackrel{1*}{\longleftrightarrow}y_h \text{ around } \Lambda\}
\, ,
$$
implies that for all $h>l$ 
\begin{align}\label{2maxdelta^2/4}
\max \left\{ 
\mu\left(  \cw{x_{h}}{y_{h}} \text{ around } \Lambda \right) , 
\mu\left(  \ccw{x_{h}}{y_{h}} \text{ around } \Lambda \right)
\right\} 
\ge 
c^2 /4
\, .
\end{align}
Additionally, considering \eqref{contra1} and 
$$
\forall h>l : ~
\{\cw{x_{h}}{y_{h}} \text{ around }\Lambda \}\subset \{\cw{x_{h}}{y_{h}} \text{ around }\Gamma \}
$$ 
yields that for all $h>l$, 
$$
\mu \left( \cw{x_{h}}{y_{h}} \text{ around }\Lambda )\right) <c^2 /4
\, .
$$
Hence, \eqref{2maxdelta^2/4} implies that for all $h>l$
$$
\mu (\ccw{x_{h}}{y_{h}} \text{ around }\Lambda)\ge c^2 /4
\, ,
$$
which, together with \eqref{h>m} and the positive association, yields 
$$
\mu\left( \ccw{x_{l+1}}{y_{l+1}} \text{ around }\Lambda, \cw{x_{l+1}}{y_{l+1}} \text{ around }\Lambda\right) 
\ge
\delta c^2/ 4
=
\delta'
\, .
$$
Given this event, Lemma \ref{cwccw} ensures the existence of a $1*$circuit around $\Lambda$, a contradiction to \eqref{0-Lambda}.
$\hfill\square$
\\

Note that the proof of this lemma relies on all five conditions of $\mu$, but, fortunately, the bounded energy condition is used only once to verify the existence of a constant $\delta>0$ such that \eqref{h>m} holds. 
Keeping this in mind will help us in the proof of Theorem \ref{mu'}, where $\mu'$ does not satisfy the bounded energy condition. 
Before we turn towards this, we obtain Theorem \ref{mu} by applying Lemma 2 to 4. 
\\


\noindent{\bf Proof of Theorem \ref{mu}: }
Let $\Delta\ssubset\Z^2 $ be large enough so that 
\begin{align}\label{0-lambda}
\mu (\Delta\stackrel{0}{\longleftrightarrow}\infty )\ge 1-c^4 /2^{5} 
\, .
\end{align}
Lemma \ref{step2} allows us to choose a square $\{ -m,\ldots ,m\}^2=:\Gamma $ with $\Delta\subset\Gamma$ such that with probability at least $c^2/2$, for any two distinct points $x,y\in\Gamma^c$, $x$ and $y$ are $1*$connected in $\Delta^c$. 
Since $\{ \cw{x_h}{y_h}\text{ around } \Delta\} \cup\{ \ccw{x_h}{y_h}\text{ around } \Delta\} = \{x_h\stackrel{1*}{\longleftrightarrow} y_h \text{ in }\Delta^c\} $, this implies that 
\begin{align}\label{maxdelta^2/4}
\max \left\{ 
\mu (\cw{x_h}{y_h}\text{ around } \Delta ),\mu (\ccw{x_h}{y_h}\text{ around } \Delta )
\right\} \ge c^2 /4 
\end{align}
for all $h\in\Z$. 
Lemma \ref{step3} thus gives the existence of nodes $x$ on the left side and $y$ on the right side of $\Gamma$ such that 
\begin{align}
\exists h>0:~ &
\mu (\cw{x_h}{y_h} \text{ around }\Delta ) 
\ge  
c^2 /4
\label{Eh>0}
\\
\exists h<0:~&
\mu (\ccw{x_h}{y_h} \text{ around }\Delta ) 
\ge  
c^2 /4 
\label{Eh<0}
\, .
\end{align}

The inequalities \eqref{maxdelta^2/4}, \eqref{Eh>0}, and \eqref{Eh<0} yield that there exists a $k\in\Z$ such that 
\begin{align}\label{hh'>0}
\mu \left( 
\cw{x_{k+1}}{y_{k+1}}\text{ around } \Delta
\right) 
, 
\mu \left( 
\ccw{x_k}{y_k}\text{ around } \Delta
\right) 
\ge 
c^2 /4 
\, .
\end{align}
Moreover, since $\{\cw{x_{k+1}}{y_{k+1}}\text{ around } \Delta\}$ and $\{\ccw{x_k}{y_k}\text{ around } \Delta\}$ are increasing events, we can conclude that
$$
\mu \left( 
\cw{x_{k+1}}{y_{k+1}}\text{ around } \Delta
, 
\ccw{x_k}{y_k}\text{ around } \Delta
\right) 
\ge 
c^4 /16 
\, .
$$ 
Thus, because of Lemma \ref{cwccw} a $1*$circuit around $\Delta$ occurs with probability at least $c^4 /16 $, a contradiction to \eqref{0-lambda} and $\mu$ cannot exist.
$\hfill\square$
\\

Only one small modification of this proof is necessary to prove the second theorem. 
\\

\noindent{\bf Proof of Theorem \ref{mu'}: } 
The strategy is to show that conditions i)-iii) and v) of Theorem \ref{mu} are satisfied and that a sufficiently close analogon to equation \eqref{h>m}, which is the only point where the bounded energy condition enters the proof, can be verified. 
 
The first three conditions of Theorem \ref{mu} are equal to the first three conditions of Theorem \ref{mu'}.

Condition v) is a consequence of conditions ii), iii), iv') and v'): 
Since the set of edges is countably infinite and the infinite boundary exists with positive probability, there exists an edge which is hit by the infinite boundary with positive probability $\eta $. 
Let $a,b$ be the nodes connected by this edge and assume without loss of generality that with probability $\eta/2$, the infinite $1*$cluster contains node $a$ and the infinite $0$cluster contains node $b$. 
Because the infinite boundary is translation invariant, shifting does not change the probability and, consequently, for all $z\in\Z^2 $ 
$$
\mu' (
z\stackrel{1*}{\longleftrightarrow}\infty 
)
\ge
\eta /2
>0
\, .
$$ 

Next, we verify a sufficiently close analogue to equation \eqref{h>m} with the notation of the proof of Lemma \ref{step3}: 

Denote by $\zeta$ the probability that an infinite $1*$cluster exists, i.e.\!,
$$
\zeta 
:=
\mu' (\Z^2\stackrel{1*}{\longleftrightarrow} \infty )
>
0
\, .
$$
Let $\Xi\ssubset\Z^2$ be large enough so that with probability at least $3\zeta /4 $, the infinite boundary exists and hits $\Xi$. 
We recall that $\Delta$ was defined in the proof of Lemma \ref{step3} as an arbitrary (but fixed) finite set of $\Z^2$. 
Take two translates $\Xi'$ and $\Xi''$ of $\Xi$ such that every node of $\Xi'$ is on the left side of $\Delta$ and every node of $\Xi''$ is on the right side of $\Delta$. 

By subadditivity of $\mu'$, the infinite boundary hits both sets $\Xi'$ and $\Xi''$ with probability at least $\zeta /4$. 
Moreover, one can find two pairs $x,x'$ and $y,y'$ of adjacent sites in $\Xi'$ resp. $\Xi''$ such that the event 
$$
\left\{ 
x\stackrel{1*}{\longleftrightarrow}\infty 
,
x'\stackrel{0}{\longleftrightarrow}\infty 
,
y\stackrel{1*}{\longleftrightarrow}\infty 
,
y'\stackrel{0}{\longleftrightarrow}\infty 
\right\} 
$$
occurs with positive probability $\epsilon $.

Take a square $[-i,i]^2 $ with $\Xi'\cup\Xi''\subset [-i,i]^2 $ such that with probability at least $\delta :=\epsilon /2$, the part of the boundary that starts between $x$ and $x'$ and ends between $y$ and $y'$ exists and is contained in $[-i,i]^2 $. 
Since the distribution of the infinite boundary is translation-invariant, for all $h\in\Z $, the event that the part of the infinite boundary starting between $x_h$ and $x'_h$ and ending between $y_h$ and $y'_h$ exists and is contained in $[-i,i]\times [-i+h,i+h] $ occurs with probability at least $\delta$, where $x_h$ is defined by $x+(0,h)$. 
Moreover, given this event, one can in fact find a $1*$path from $x_h$ to $y_h$ in $[-i,i]\times [-i+h,i+h] $. 
This, obviously, implies that for all $h>2 \max \{i,m\} $ 
$$
\mu \left( \cw{x_h }{y_h } \text{ around } \Delta \right) 
\ge 
\delta
\, ,
$$
which is sufficiently close to \eqref{h>m}. 
$\hfill\square$ 
\\

For the sake of completeness we give the proof of the first corollary:
\\

\noindent{\bf Proof of Corollary \ref{rho}: } 
Because of the Burton-Keane uniqueness theorem at most one infinite $1*$cluster as well as at most one infinite $0$cluster exist. 
We assume for contradiction that both of them coexist with strictly positive probability. 
The ergodicity yields that this event occurs with $\rho$-probability one. 
So, all conditions of Theorem \ref{mu'} are satisfied and the contradiction is shown. 
$\hfill\square$ 
\\

\section{A single infinite $1*$cluster has unbounded width}

From now on let $\nu$ be a probability measure on $\{ 0,1\}^{\Z^2}$ satisfying the bounded energy condition. 

Our aim is to show that given the existence of an infinite $1*$cluster, one can find infinitely many infinite self-avoiding $1*$paths. 
To this end, we first have to check the measurability of the latter event
, where the corresponding $\sigma$-algebra is generated by the cylinder sets.

\begin{lem}\label{lem:tailevent}
The number $A$ of infinite self-avoiding $1*$paths is tail measurable. 
\end{lem}

\pf
The statement is a direct consequence of the identity 
$$
\{ A\ge n \}
=
\bigcup_{l\in\N} \bigcap_{k\ge l }
\bigcup_{m\ge k }\bigcap_{i\ge m } \{ A_{k,i} \ge n \}
\, ,
$$
which holds for all $n\in\N$. 
Here, $A_{k,i}$ is the maximal number of disjoint $1*$paths in $\{ -i,\ldots ,i\}^2\setminus\{ -k,\ldots ,k\}^2$ from $\partial^* \{ -k,\ldots ,k\}^2$ to $\partial^* \{-i+1,\ldots ,i-1\}^2 $. 
$\hfill\square$
\\

Next, we show that configurations with a given number of disjoint infinite self-avoiding $1*$paths exhibit a so-called necklet with this number of $1$pearls around any finite set. 
\begin{defi}
Let $N\in\N$, $\sigma\in\{ 0,1\}^{\Z^2}$, and $\Gamma\ssubset\Z^2$.
We call $C$ a necklet with $N$ $1$pearls around $\Gamma$ with respect to $\sigma $ if $C$ is a circuit around $\Gamma $  with $|C\cap\sigma^{-1}(1) |=N$. 
\end{defi}
Note that the proof of the following existence statement is more or less a direct consequence of the well-known max-flow min-cut theorem of Ford and Fulkerson; cf \cite{FF}.
Since this is the only point where the max-flow min-cut theorem (and its notation) is needed, we use the original notation of \cite{FF} without defining it. 

\begin{lem}[Bottleneck lemma]\label{lem:necklet}
Let $\sigma\in\{ 0,1\}^{\Z^2} $ be a configuration that possesses exactly $N$ disjoint infinite self-avoiding $1*$paths. 
Then, for all $\Gamma\ssubset\Z^2$, there exists a necklet with $N$ $1$pearls around $\Gamma$. 
\end{lem}

\pf 
Fix an arbitrary $\Gamma\ssubset\Z^2$ and a configuration $\sigma$ such that one can find exactly $N$ disjoint infinite self-avoiding $1*$paths with respect to $\sigma$. 
Let the set $S$ of sources be the $*$boundary of a square $\{-s ,\ldots ,s\}^2$ large enough so that it contains $\Gamma$ and $N$ disjoint infinite self-avoiding $1*$paths start in this square. 
Furthermore, the set $T$ of sinks is defined as the $*$boundary of a square $\{-t ,\ldots ,t\}^2$ large enough so that $S\subset \Int T $ and there exist $N$ disjoint self-avoiding $1*$paths from $S$ to $T$. 
The set of intermediate nodes $R$ is $\Int T\setminus (S\cup\Int S)$. 
An undirected arc $\{ x,y\}$ connects $x$ and $y$ if and only if these two nodes belong to $ R \cup S\cup T$ and are $*$adjacent. 
%
We define the capacity function $c(.,.)$ of an arc $\{ x,y\}$ as
$$
c(x,y)= 
\begin{cases}
1 & \text{ if }x,y\subset\sigma^{-1}(1) , \\
0 & \text{ otherwise. }
\end{cases}
$$
Consequently, since there are $N$ disjoint self-avoiding $1*$paths from $S$ to $T$, the maximal flow value of this network is $N$. 
Applying the max-flow min-cut theorem, see \cite[page 11 plus section 7 and 10]{FF}, shows the existence of a cut $C$ separating $S$ from $T$, whose cut capacity is $N$. 

Let $B$ be the set of nodes that are connected to $S$ by a $*$path not using an arc of the cut $C$. 
Because the capacity of $C$ is $N$ and $N$ disjoint self-avoiding $1*$paths connect $S$ to $T$, there exists only one cluster $B^o$ in $\partial^* B$ that is also a self-avoiding circuit around $S$. 
By definition of $B$, the set $B^i$ of all nodes in the interior of $B^o$ $*$adjacent to $B^o$ form a circuit around $\Int S$. 
Furthermore, a node of $B^i$ and a node of $B^o$ $*$adjacent to each other are also connected by an arc of the cut $C$. 

Now we are ready to construct the necklet: 
First, take the set $D$ of nodes of $B^i\cap\sigma^{-1}(0)$ and combine it with the set $E$ of nodes in $B^i\cap\sigma^{-1}(1)$ $*$adjacent to $B^o\cap\sigma^{-1}(1)$. 
Since $N$ disjoint self-avoiding $1*$paths connect $S$ to $T$ and the capacity of $C$ is $N$, the set $E$ consists of exactly $N$ nodes. 
By definition, the set $F$ of nodes in $B^o$ $*$adjacent to $B^i\setminus (D\cup E)$ is contained in $\sigma^{-1}(0)$. 
A moment's thought reveals that $D\cup E\cup F$ is a necklet with $N$ $1$pearls around $\Gamma$ with respect to $\sigma $.
$\hfill\square$
\\

Now, let us gain some insight into the structure of infinite $1*$clusters under fairly general conditions on the measure. 
\\

\noindent{\bf Proof of Theorem \ref{nu}: }
Since $\left( \{0,1 \} ,\Pot (\{0,1 \})\right)$ is a perfect space, Theorem 3.3 of \cite{S} implies that $\nu $ is a Gibbs measure for a suitable specification $(\gamma_ \Lambda)_{\Lambda\ssubset\Z^2 }$. 
Since $\nu$ satisfies the bounded energy condition, there exist constants $c_n>0$ such that 
$
\gamma_\Lambda (\eta |\xi )\ge c_n
$ 
for $\nu$-almost all configurations $\xi\in\{ 0,1\}^{\Z^2}$, whenever $|\Lambda|=n$ and $\eta$ is a local configuration on $\Lambda$. 
Applying the extreme decomposition theorem (7.26) of \cite{Geo} yields that the bounded energy condition holds for $\pr_{\nu }$-almost all extreme Gibbs measure specified by $(\gamma_ \Lambda)_{\Lambda\ssubset\Z^2 }$.
So, we may assume without loss of generality that $\nu $ is trivial on the tail $\sigma$-field. 

We further assume without loss of generality $\nu ( \Z^2  \stackrel{1*}{\longleftrightarrow} \infty )>0$. 
The triviality of $\nu$ on the tail $\sigma$-field then implies $\nu ( \Z^2  \stackrel{1*}{\longleftrightarrow} \infty )=1$. 
Consequently, we just have to verify that infinitely many disjoint infinite self-avoiding $1*$paths $\nu $-almost surely exist.  
The proof of the other statement is similar. 

By assumption, the number $A$ of infinite self-avoiding $1*$paths $\nu$-almost surely is at least one. 
We will show that $\nu(A=\infty)=1$ or, equivalently, that $\nu (A=N ) =0$ for all $N\ge 1$. 

Suppose the contrary. 
Tail triviality, together with Lemma \ref{lem:tailevent}, implies the existence of some $N\ge 1$ with $\nu (A=N ) =1$. 
Because $\nu$ satisfies the bounded energy condition we can choose an $\epsilon >0 $ such that 
\begin{align}\label{Sle5N}
\nu \left(  
\eta \text{ on } S
\bigg| 
\xi \text{ off } S
\right) 
\ge
\epsilon
\end{align}
for all  $S\subset\Z^2$ with $|S|\le 5N$, $\eta\in \{ 0\}^S $ and for $\nu $-a.e.\! $\xi\in\{ 0,1\}^{S^c }$. 
Let $\Gamma \ssubset\Z^2$ be large enough so that 
\begin{align}\label{DToInfty}
\nu (\Gamma\stackrel{1*}{\longleftrightarrow}\infty) >1-\epsilon /4
\, .
\end{align}
The bottleneck lemma ensures the $\nu$-almost-sure existence of a necklet with $N$ $1$pearls around $\Gamma $. 
Let $\Delta\ssubset\Z^2$ be large enough so that with probability at least $1-\epsilon /2$, there exists a $0$necklet with $N$ $1$pearls around $\Gamma\cup\partial^*\Gamma $ in $\Delta$. 

Denote by $C$ the maximal $0$necklet with $N$ $1$pearls around $(\Gamma\cup\partial^*\Gamma )$ in $\Delta $; 
if it does not exist $C$ is $\emptyset$. 
Hence, $\Int C$ is a well-defined random set, which is determined from outside. 
Let $S$ be the set of nodes in $\Int C$ $*$adjacent to a $1$pearl of $C$, under the condition $C\neq\emptyset $. 
Otherwise $S$ is $\emptyset$.  
Once again $S$ is a well-defined random set, which is determined from outside of $\Int C$ and $|S|\le 5N$ always holds. 
If $C\neq\emptyset$ and all spins of $S$ take the value zero, a $0$circuit around $\Gamma $ exists. 
Hence, the inequality \eqref{Sle5N} yields that the existence of a $0$circuit around $\Gamma $ in $\Delta$ has probability at least $(1-\epsilon/2) \epsilon$, a contradiction to \eqref{DToInfty}. 
Consequently, $\nu (A\in\N ) =0$.  
$\hfill \square$
\\

\textbf{Acknowledgments:} 
It is my pleasure to thank H.-O. Georgii for stimulating discussions and helpful remarks. 
In particular H.-O. Georgii pointed out an error in an earlier version of the proof of Theorem \ref{nu}. 
Also I would like to thank an anonymous referee for hinting at the papers of Sheffield \cite{Sheff} and Bollobas and Riordan \cite{BR}. 
Last but not least I would like to cordially thank Ursula and J\" urgen Carstens for funding and never-ending support.

\newpage 

\renewcommand{\thesection}{}

\setlength{\parindent}{0cm}

\end{document}